\documentclass[12pt,twoside]{article}
\title{The
projectively resolving of some classes  over a direct product of
rings}
\date{}
\usepackage{amsfonts}
\usepackage{amsmath}
\usepackage{amssymb}
\usepackage{latexsym}
\usepackage{graphicx}
\newtheorem{thm}{\bf Theorem}[section]

\newtheorem{lem}[thm]{\bf Lemma}
\newtheorem{prop}[thm]{\bf Proposition}
\newtheorem{defn}[thm]{\bf Definition}

\newtheorem{exmp}[thm]{\bf Example}

\catcode`\ç=13
\defç{\c{c}}
\catcode`\é=13
\defé{\'e}
\catcode`\à=13
\defà{\`a}
\catcode`\è=13
\defè{\`e}
\catcode`\â=13
\defâ{\^a}
\catcode`\ù=13
\defù{\`u}
\catcode`\ê=13
\defê{\^e}
\catcode`\î=13
\defî{\^\i}
\catcode`\ô=13
\defô{\^o}

\def\proof{{\parindent0pt {\bf Proof.\ }}}

\def\Im{{\rm Im}}

\begin{document}
\thispagestyle{empty}
\maketitle \vspace*{-2cm}
\begin{center}{\large\bf Najib Mahdou and Mohammed Tamekkante}

\bigskip
\small{Department of Mathematics, Faculty of Science and Technology
of Fez,\\ Box 2202, University S. M.
Ben Abdellah Fez, Morocco, \\ mahdou@hotmail.com  \\
tamekkante@yahoo.fr}
\end{center}

\bigskip\bigskip
\noindent{\large\bf Abstract.} In this paper, we study the resolving
of  $\mathcal{SGP}(-)$ and $\mathcal{SGF}(-)$, the classes of all
strongly Gorenstein projective and flat modules respectively, over a
 direct product of commutative rings.
\bigskip

\small{\noindent{\bf Key Words.} Classical homological dimensions of
modules; Gorenstein homological dimensions of modules and of rings;
resolving classes; (strongly)Gorenstein projective, injective, and
flat modules.}

\begin{section}{Introduction}
Throughout this paper, all rings are commutative with identity
element, and all modules are unital.\\

$\mathbf{\quad Setup\; and\; Notation:}$ Let $R$ be a ring, and let
$M$ be an $R$-module. We use $pd_R(M)$, $id_R(M)$ and $fd_R(M)$ to
denote, respectively, the classical projective, injective and flat
dimensions of $M$. By $gldim(R)$ and $wdim(R)$ we denote,
respectively, the classical global and weak dimensions
of R.\\

 Recently in \cite{Bennis and Mahdou2}, the authors started the study of
global Gorenstein dimensions of rings, which are called, for a
commutative ring $R$, projective, injective, and weak dimensions of
$R$; denoted by $GPD(R)$, $GID(R)$, and $G - wdim(R)$, respectively;
and, respectively, defined as follows:\bigskip

$\begin{array}{cccc}
  1) & GPD(R) & = & sup\{ Gpd_R(M)\mid M$ $R-module\} \\
  2) & GID(R) & = & sup\{ Gid_R(M)\mid M$ $R-module\} \\
  3) & G-wdim(R) & = & sup\{ Gfd_R(M)\mid M$ $R-module\}
\end{array}$
\\

They proved that, for any ring R, $ G-wdim(R)\leq GID(R) = GPD(R)$
(\cite[Theorems 2.1 and 2.11]{Bennis and Mahdou2}). So, according to
the terminology of the classical theory of homological dimensions of
rings, the common value of $GPD(R)$ and $GID(R)$ is called
Gorenstein global dimension of $R$, and denoted by $G-gldim(R)$.\\
They also proved that the Gorenstein global and weak dimensions are
refinement of the classical global  and weak dimensions of rings.
That is : $G-gldim(R) \leq gldim(R)$ and $G-wdim(R)\leq wdim(R)$
with equality if $wdim(R)$ is finite (\cite[Proposition 2.12]{Bennis
and Mahdou2}).\bigskip

$\mathbf{Notation.}$ By  $\mathcal{P}(R)$, $\mathcal{I}(R)$ and
$\mathcal{F}(R)$ we denote the classes of all projective, injective
and fat $R$-modules respectively; and by
 $\mathcal{GP}(R)$,
$\mathcal{GI}(R)$ and $\mathcal{GF}(R)$, we denote the classes of
all strongly Gorenstein projective, injective and fat $R$-modules
respectively. Furthermore, we let $\mathcal{SGP}(R)$ and
$\mathcal{SGF}(R)$
denote the classes of all  Gorenstein  projective and flat $R$-modules, respectively.\\

In  \cite{Holm}, Holm introduced the notion of resolving classes as
follows:
\begin{defn}(\cite[Definition 1.1]{Holm}) For any class $\mathcal{X}$ of $R$-modules.
\begin{description}
    \item[a/]We call $\mathcal{X}$ projectively resolving if $\mathcal{P}(R)\subseteq \mathcal{X}$, and for every short exact sequence
$0 \longrightarrow X' \longrightarrow X \longrightarrow X"
\longrightarrow 0$ with $X" \in \mathcal{X}$ the conditions $X' \in
\mathcal{X}$ and $X \in \mathcal{X}$ are equivalent.
    \item[b/]We call $\mathcal{X}$ injectively resolving if $\mathcal{I}(R)\subseteq \mathcal{X}$, and for every short exact sequence
$0 \longrightarrow X' \longrightarrow X \longrightarrow X"
\longrightarrow 0$ with $X' \in \mathcal{X}$ the conditions $X" \in
\mathcal{X}$ and $X \in \mathcal{X}$ are equivalent.
\end{description}
\end{defn}
In \cite{Holm} again, Holm prove that  the class $\mathcal{GP(R)}$
is projectively resolving and closed under arbitrary direct sums and
under direct summands (\cite[Theorem 2.5]{Holm}); and dually, the
class $\mathcal{GI(R)}$ is injectively  resolving and closed under
arbitrary direct products and under direct summands (\cite[Theorem
2.6]{Holm}). He also prove that, if $R$ is coherent, then the class
$\mathcal{GF(R)}$ is projectively resolving and closed under direct
summands (\cite[Theorem 3.7]{Holm}). \bigskip

In section 2, we are mainly interested to study the resolving of the
class $\mathcal{SGP}(-)$ over a direct
product of rings.\\

Similarly, in section 3, we study the resolving of the class
$\mathcal{SGP}(-)$  and over a direct product of rings. We find
after, some connections between the projectively resolving of the
classes $\mathcal{SGP}(-)$ and $\mathcal{SGF}(-)$ under some
conditions.\\

\end{section}
\begin{section}{The resolving of the class $\mathcal{SGP}(-)$ over a 
direct product of rings}

The aim of this section is to study the connection between the
projectively resolving of the class $\mathcal{SGP}(-)$ over a
 family of commutative rings and direct product of this family.

  The main result in this section is:
\begin{thm}\label{SGP product} Let $\{R_i\}_{i=1,..,m}$ be a family of rings with finite Gorenstein global dimensions.
Then,
   $\mathcal{SGP}(\displaystyle \prod_{i=1}^mR_i)$ is projectively resolving if, and only
 if,
 $\mathcal{SGP}(R_i)$ is projectively resolving for each
 $i=1,..,m$.
 \end{thm}
To prove this theorem, we need some results.\\
First, we need the notion of strongly Gorenstein projective module,
which is introduced in \cite{Bennis and Mahdou1}, in order to
characterize Gorenstein projective modules.

\begin{defn}
 A module $M$ is said to be strongly Gorenstein projective ($SG$-projective for short), if
there exists an exact sequence  of the form:
$$\mathbf{P}=\ \cdots\rightarrow P\stackrel{f}\rightarrow
P\stackrel{f} \rightarrow P\stackrel{f}
         \rightarrow P \rightarrow\cdots$$ where $P$ is a projective $R$-module and $f$ is an endomorphism of $P$, such
          that  $M \cong \Im(f)$ and such that $\mathbf{Hom}(-,Q)$ leaves the sequence $\mathbf{P}$ exact whenever
           $Q$ is a projective module.
\end{defn}

These strongly Gorenstein projective modules have a simple
characterization, and they are used to characterize Gorenstein
projective modules. We have:

\begin{prop}\cite[Propositions 2.9]{Bennis and Mahdou1}
  A module M is strongly Gorenstein
projective if, and only if, there exists a short exact sequence of
modules: $$0\longrightarrow M \longrightarrow P \longrightarrow M
\longrightarrow 0$$ where $P$ is projective and $Ext(M,Q) = 0$ for
any projective module $Q$.
\end{prop}

\begin{lem}\label{lem SGP}
Let $R$ be a  ring.
 The class $\mathcal{SGP}(R)$ is $\mathcal{P}(R)$-resolving
 if, and only if, every Gorenstein projective $R$-module is strongly
    Gorenstein projective.
\end{lem}
\proof First, assume that  $\mathcal{SGP}(R)$ is
$\mathcal{P}(R)$-resolving. By \cite[Proposition 2.2]{Bennis and
Mahdou1}, $\mathcal{SGP}(R)$ is stable by direct sums. Then, using
 \cite[Proposition 1.4]{Holm},
$\mathcal{SGP}(R)$ is closed under direct summands. On the other
hand, from \cite[Theorem 2.7]{Bennis and Mahdou1}, every
$G$-projective $R$-module is direct summand of an $SG$-projective
$R$-module. Consequently,
 every
$G$-projective $R$-module is $SG$-projective. Conversely, assume
that every $G$-projective $R$-module is $SG$-projective. Then we
have $\mathcal{GP}(R) \subseteq \mathcal{SGP}(R) \subseteq
\mathcal{GP}(R)$. So, $\mathcal{GP}(R)= \mathcal{SGP}(R)$.
Consequently, from \cite[Theorem 2.5]{Holm}, the class
$\mathcal{SGP}(R)$ is projectively resolving.$\blacksquare$
\begin{thm}{\cite[Theorem 2.7]{Bennis and Mahdou1}} A  module is Gorenstein projective
if, and only if, it is a direct summand of a strongly Gorenstein
projective  module.
\end{thm}

$\mathbf{Proof \; of \; Theorem \; \ref{SGP product}}.$
By induction on m, we may assume $m = 2$.\\
Assume that  $\mathcal{SGP}(R_1\times R_2)$ is projectively
resolving. By Lemma \ref{lem SGP}, to prove that
$\mathcal{SGP}(R_1)$ is projectively resolving, it remains to prove
that every $G$-projective $R_1$-module is an $SG$-projective
$R_1$-module. So, let $M$ be a $G$-projective $R_1$- module.
$M\times 0$ is  an $R_1\times R_2$-module (see \cite[Page
101]{Berrick}). First, we claim that  $M\times 0$ is a
$G$-projective  $R_1\times R_2$-module. The $R_1$-module $M$ is a
direct summand of an $SG$-projective $R_1$-module $N$ (\cite[Theorem
2.7]{Bennis and Mahdou1}). For such module, there is a short exact
sequence of $R_1$-modules:
$$0\longrightarrow N\longrightarrow P \longrightarrow N \longrightarrow 0$$  where $P$ is a projective $R_1$-module.
Hence, we have a  short sequence of $R_1\times R_2$-modules:
$$(\ast)\quad 0\longrightarrow N\times 0\longrightarrow P\times 0 \longrightarrow
N\times 0 \longrightarrow 0$$ and  $P\times 0$ is a projective
$R_1\times R_2$-module. But $G-gldim(R_1\times R_2)$ is finite.
Then, there is an integer $i>0$ such that $Ext_{R_1\times R_2}^i(
N\times 0,Q)=0$ for each projective $R_1\times R_2$-module $Q$
(since $Gpd_{R_1\times R_2}(N\times 0)<\infty$ and by \cite[Theorem
2.20]{Holm}). From $(\ast)$ we deduce that $Ext_{R_1\times R_2}(
N\times 0,Q)=0$. Then, by \cite[Proposition 2.9]{Bennis and
Mahdou1}, $N\times0$ is an $SG$-projective $R_1\times R_2$-module.
So, $M\times0$ is a $G$-projective $R_1\times R_2$-module (since it
is a direct summand of $N\times0$ as $R_1\times R_2$-modules and by
\cite[Theorem 2.5]{Holm}). Now, we claim that $M$ is an
$SG$-projective $R_1$-module. From Lemma \ref{lem SGP}, the
$R_1\times R_2$-module $M\times 0$ is $SG$-projective (since
$\mathcal{SGP}(R_1\times R_2)$ is
    resolving and $M\times0$ is a $G$-projective
$R_1\times R_2$-module). Thus,  there exists
    a short exact sequence of $R_1\times R_2$-modules: $$(\star)\quad 0
    \longrightarrow M\times 0\longrightarrow P \longrightarrow
    M \times 0 \longrightarrow 0$$ where $P$ is a  projective $R_1\times
    R_2$-module. Now, we tensor $(\star)$ by $-\otimes R_1$ and we obtain the short exact sequence of $R_1$-modules (see that $R_1$ is a projective
$R_1\times R_2$-module):$$(\star\star)\quad 0
    \longrightarrow M\times 0 \displaystyle\otimes _{R_1\times R_2}R_1 \longrightarrow P\displaystyle\otimes _{R_1\times R_2}R_1 \longrightarrow
    M\times 0 \displaystyle\otimes _{R_1\times R_2}R_1 \longrightarrow
    0$$ But  $M\times 0 \displaystyle\otimes _{R_1\times
    R_2}R_1\cong M\times0\displaystyle\otimes_{R_1\times
    R_2}(R_1\times R_2)/(0\times R_2) \cong M\times 0\displaystyle\cong_RM$ (isomorphism of $R$-modules). Then, we can
    write $(\star\star)$ as :$$ 0 \longrightarrow M
    \longrightarrow P\displaystyle\otimes _{R_1\times
    R_2}R_1 \longrightarrow M \longrightarrow 0$$ It is clear that $P\displaystyle\otimes
    _{R_1\times
    R_2}R_1$ is a projective $R_1$-module. Furthermore, by \cite[Theorem 2.20]{Holm}, $Ext_{R_1}(M,F)=0$ for every $R_1$-module projective
      $F$ since $M$ is a $G$-projective $R_1$-module. So, by \cite[Proposition 2.9]{Bennis and Mahdou1}, $M$ is an
    $SG$-projective $R_1$-module, as desired.\\
   Similarly, we can prove that $\mathcal{SGP}(R_2)$ is
   projectively  resolving.\\

 Conversely, assume that $\mathcal{SGP}(R_i)$
for $i=1,2$ are projectively resolving and let
    $M$ be a $G$-projective $R_1\times R_2$-module. We claim that $M$ is an $SG$-projective $R_1\times R_2$-module. We have the isomorphism of $R_1\times R_2$-modules: $$M\cong M\otimes_{R_1\times
    R_2}R_1\times R_2\cong M\otimes_{R_1\times
    R_2}(R_1\times0\oplus 0\times R_2)\cong M_1\times M_2$$ where
    $M_i=M\otimes_{R_1\times
    R_2} R_i$ for $i=1,2$ (for more details see \cite[p.102]{Berrick}). By \cite[Lemma 3.2]{Bennis and Mahdou3}, for each $i=1,2$,  $M_i$ is
    a $G$-projective $R_i$-module. Then, by Lemma \ref{lem SGP}, $M_i$ is an $SG$-projective $R_i$-module since $\mathcal{SGP}(R_i)$ is projectively resolving for each $i=1,2$.
     On the other hand, the family $\{R_i\}_{i=1,2}$ of rings satisfies the conditions of
    \cite[Lemma 3.3]{Bennis and Mahdou3} (by \cite[Corollary 2.10]{Bennis and Mahdou2}  since $G-gldim(R_i)$ is finite  for each $i=1,2$). Thus,   $M=M_1\times M_2$ is
   an  $SG$-projective  $R_1\times R_2$-module and the desired result follows from  Lemma \ref{lem SGP}.$\blacksquare$\

\begin{prop}
Let $R$ be a ring such that $G-gldim(R)\leq 1$ and
$\mathcal{SGP}(R)$ is projectively resolving. Then, $R$ is coherent.
\end{prop}
\proof Let $R$ be a ring  which satisfies the  conditions above and
let $I$ be a finitely generated ideal of $R$. From the short exact
sequence $0 \longrightarrow I \longrightarrow R \longrightarrow
R/I\longrightarrow 0$ we see that $I$ is $G$-projective since
$G-gldim(R)\leq 1$. Then, by Lemma \ref{lem SGP}, $I$ is
$SG$-projective since $\mathcal{SGP}(R)$ is projectively resolving.
Hence, from \cite[Proposition 3.9]{Bennis and Mahdou1}, $I$ is
finitely
presented and so $R$ is coherent.$\blacksquare$\\

Now, we construct a  non-Noetherian rings $\{R_i\}_{i\geq 1}$ such
that $\mathcal{SGP}(R_i)$ is projectively resolving, $G-gldim(R_i) =
i$ and $wdim(R_i) = \infty$ for all $i \geq 1$, as follows:
\begin{exmp} Consider a non-semisimple quasi-Frobenius ring
$R=K[X]/(X^2)$ where $K$ is a field, and a non- Noetherian
hereditary ring $S$. Then, for every positive integer $n$, we have:
$\mathcal{SGP}(R \times S[X_1,X_2, ...,X_n])$ is projectively
resolving, $G-gldim(R \times S[X_1,X_2, ...,X_n]) = n + 1$ and
$wdim(R \times S[X_1,X_2, ...,X_n]) = \infty$.

\end{exmp}
\proof From \cite[Example 3.4]{Bennis and Mahdou3}, only the
projectively  resolving of $R \times S[X_1,X_2, ...,X_n]$ need an
argument. It is clear that $\mathcal{SGP}(S[X_1,X_2,...,X_n])$ is
projectively resolving since $wdim(S[X_1,X_2,...,X_n])$ is finite
 (by the Hilbert Syzygies's Theorem) and since, by \cite[Proposition 2.27]{Holm}, $\mathcal{SGP}(S[X_1,X_2,...,X_n])=\mathcal{P}(S[X_1,X_2,...,X_n])$ (the class of all projective modules). On the other hand, from \cite[Corollary 3.9]{Ouarghi}, $\mathcal{SGP}(R)$ is
projectively resolving. Thus, from Theorem \ref{SGP product},
$\mathcal{SGP}(R \times S[X_1,X_2, ...,X_n])$ is projectively
resolving.$\blacksquare$

\end{section}
\begin{section}{The resolving of the class $\mathcal{SGF}(-)$ over a direct product of rings}

The aim of this section is to study the connection between the
projectively resolving of the class $\mathcal{SGF}(-)$ over a family
of rings and over a direct product of this family.\\

 Our main result in this section is:

\begin{thm}\label{SGF product}
Let $\{R_i\}_{i=1,..,m}$ be a family of coherent rings with finite
Gorenstein weak dimensions. Then,
 $\mathcal{SGF}(\displaystyle \prod_{i=1}^mR_i)$ is projectively resolving if, and only
 if,
 $\mathcal{SGF}(R_i)$ is projectively resolving for each $i=1,..,m$.
 \end{thm}

To prove this theorem, we need some results.\\
First, we need the notion of strongly Gorenstein flat modules, which
is introduced in \cite{Bennis and Mahdou1} to characterize the
Gorenstein projective modules.
\begin{defn}
A module $M$ is said to be strongly Gorenstein flat ($SG$-flat for
short), if there exists an exact sequence  of the form:
$$\mathbf{F}=\ \cdots\rightarrow F\stackrel{f}\rightarrow F \stackrel{f}\rightarrow
F \stackrel{f}\rightarrow F \rightarrow\cdots$$ where $F$ is flat
and $f$ is an endomorphism of $F$ such that $M \cong \Im(f)$ and
such that $-\bigotimes I$ leaves the sequence $\mathbf{F}$ exact
whenever $I$ is an injective module.
\end{defn}
These strongly Gorenstein flat modules have a simple
characterization, and they are used to characterize the Gorenstein
projective modules. We have:
\begin{prop}\cite[Proposition 3.6]{Bennis and Mahdou1}
A module M is strongly Gorenstein flat  if, and only if, there
exists a short exact sequence of modules: $0\longrightarrow M
\longrightarrow F \longrightarrow M \longrightarrow 0$ where $F$ is
flat and $Tor(M,I) = 0$ for any injective  module $I$.
\end{prop}

 \begin{lem}\label{lem SGF}
Let $R$ be a coherent ring. The class $\mathcal{SGF}(R)$ is
 projectively resolving if, and only if, every Gorenstein flat
$R$-module is strongly Gorenstein flat.
\end{lem}

\proof First, assume that  $\mathcal{SGF}(R)$ is projectively
resolving. By \cite[Proposition 3.4]{Bennis and Mahdou1},
$\mathcal{SGF}(R)$ is stable by direct sums. Then, using
 \cite[Proposition 1.4]{Holm},
$\mathcal{SGF}(R)$ is closed under direct summands. On the other
hand, from \cite[Theorem 3.5]{Bennis and Mahdou1}, every $G$-flat
$R$-module is a direct summand of an $SG$-flat $R$-module.
Consequently,
 every
$G$-flat $R$-module is $SG$-flat. Conversely, assume that every
$G$-flat $R$-module is $SG$-flat. Then we have $\mathcal{GF}(R)
\subseteq \mathcal{SGP}(R) \subseteq \mathcal{GP}(R)$. So,
$\mathcal{GF}(R)= \mathcal{SGF}(R)$. Consequently, from
\cite[Theorem 3.7]{Holm}, the class $\mathcal{SGF}(R)$ is
projectively resolving (since $R$ is coherent).$\blacksquare$\\

$\mathbf{Proof \; of \; Theorem \; \ref{SGF product}.}$
By induction on m, we may assume $m = 2$.\\
First note that $R_1\times R_2$ since $R_1$ and $R_2$ are coherents.\\

 Assume that $\mathcal{SGF}(R_1\times R_2)$
is projectively resolving. By Lemma \ref{lem SGF}, to prove that
$\mathcal{SGP}(R_1)$ is projectively resolving remains to prove that
every $G$-flat $R_1$-module is an $SG$-flat $R_1$-module (since
$R_1$ is coherent). So, let $M$ be a $G$-flat $R_1$-module. We claim
that $M$ is an $SG$-flat $R_1$-module.   $M\times 0$ is an
$R_1\times R_2$-module (see \cite[Page 101]{Berrick}). First, we
claim that $M\times 0$ is a $G$-flat $R_1\times R_2$-module. The
$R_1$-module $M$ is a direct summand of an $SG$-flat $R_1$-module
$N$ (\cite[Theorem 3.5.]{Bennis and Mahdou2}). For such module,
there is a short exact sequence of $R_1$-modules:
$$0\longrightarrow N\longrightarrow F \longrightarrow N \longrightarrow 0$$  where  $F$ is a flat $R_1$-module.
Hence, we have a  short sequence of $R_1\times R_2$-modules
$$(\ast)\quad 0\longrightarrow N\times 0\longrightarrow F\times 0 \longrightarrow
N\times 0 \longrightarrow 0$$ and  $F\times 0$ is a flat $R_1\times
R_2$-module (\cite[Lemma 3.7]{Bennis and Mahdou3}). But
$G-wdim(R_1\times R_2)$ is finite. Then, there is an integer $i>0$
such that $Tor_{R_1\times R_2}^i( N\times 0,I)=0$ for each injective
$R_1\times R_2$-module $I$. From $(\ast)$ we deduce that
$Tor_{R_1\times R_2}( N\times 0,I)=0$. Then, by \cite[Proposition
3.6]{Bennis and Mahdou1}, $N\times0$ is an $SG$-flat $R_1\times
R_2$-module. So, from \cite[Theorem 3.7]{Holm}, $M\times0$ is a
$G$-flat $R_1\times R_2$-module (since it is a direct summand of
$N\times0$ as $R_1\times R_2$-modules and since $R_1\times R_2$ is
coherent). Now, we claim that $M$ is an $SG$-flat $R_1$-module. The
$R_1\times R_2$-module $M\times 0$ is $SG$-flat (since
$\mathcal{SGF}(R_1\times R_2)$ is
    resolving and by Lemma \ref{lem SGF}). Then,  there exists
    a short exact sequence of $R_1\times R_2$-modules: $$(\star)\quad 0
    \longrightarrow M\times 0\longrightarrow F \longrightarrow
    M \times 0 \longrightarrow 0$$ where $F$ is a flat $R_1\times
    R_2$-module. Now, we tensor $(\star)$ by $-\otimes_{R_1\times R_2} R_1$ (see that $R_1$ is a projective $R_1\times R_2$-module), we obtain the short exact sequence of $R_1$-modules:
 $$(\star\star)\quad 0
    \longrightarrow M\times 0 \displaystyle\otimes _{R_1\times R_2}R_1 \longrightarrow F\displaystyle\otimes _{R_1\times R_2}R_1 \longrightarrow
    M\times 0 \displaystyle\otimes _{R_1\times R_2}R_1 \longrightarrow
    0$$ But  $M\times 0 \displaystyle\otimes _{R_1\times
    R_2}R_1\cong M\times0\displaystyle\otimes_{R_1\times
    R_2}(R_1\times R_2)/(0\times R_2) \cong M\times 0\displaystyle\cong_RM$ (isomorphism of $R$-modules). Then, we can
    write $(\ast\ast)$ as :$$ 0 \longrightarrow M
    \longrightarrow F\displaystyle\otimes _{R_1\times
    R_2}R_1 \longrightarrow M \longrightarrow 0$$ It is clear that $F\displaystyle\otimes
    _{R_1\times
    R_2}R_1$ is a flat $R_1$-module. So, $M$ is an
    $SG$-flat $R_1$-module (since $G-wdim(R_1)$ is finite and by the same argument as above), as
    desired.\\
   Similarly, we can prove the projectively resolving of the class $\mathcal{SGP}(R_2)$.\\

 Conversely, assume that $\mathcal{SGP}(R_i)$
for, $i=1,2$, are projectively resolving and let
    $M$ be a $G$-flat $R_1\times R_2$-module. We have the isomorphism of $R_1\times R_2$-modules:
$$M\cong M\otimes_{R_1\times
    R_2}R_1\times R_2\cong M\otimes_{R_1\times
    R_2}(R_1\times0\oplus 0\times R_2)\cong M_1\times M_2$$ where
    $M_i=M\otimes_{R_1\times
    R_2} R_i$ for $i=1,2$ (see \cite[Page 102]{Berrick}). By \cite[Proposition 3.10]{Holm}, for each $i=1,2$,  $M_i$
    is a
    $G$-flat $R_i$-module. Then, $M_i$ is an $SG$-flat $R_i$-module  (since $\mathcal{SGF}(R_i)$ is projectively  resolving
for $i=1, 2$).\\
    Now, let $I$ be an injective $R_1$-module and put $n=G-wdim(R_1)$. Then, for
    every $R_1$-module $K$ we have $Gfd_R(K)\leq n$. Therefore, $Tor_{R_1}^{n+1}(K,I)=0$. Thus, $fd_{R_1}(I)\leq
    n$. Similarly, we can prove that every injective $R_2$-module has a finite flat
    dimension. Thus, the family $\{R_i\}_{i=1,2}$ of rings satisfies the conditions of
    \cite[Lemma 3.6]{Bennis and Mahdou3}. Hence,   $M=M_1\times M_2$ is
   an  $SG$-flat  $R_1\times R_2$-module and the desired result follows from  Lemma \ref{lem SGF} (since $R_1\times R_2$ is coherent).$\blacksquare$\\

   Recall that a ring $R$ is called perfect if every flat $R$-module is
   projective \cite{Bass}.\\

   The next result gives us a connection between the projectively
   resolving of the classes $\mathcal{SGP}(-)$ and $\mathcal{SGF}(-)$ over the some
   ring under some conditions.\\

   \begin{prop}
   Let $R$ be a coherent ring with  finite Gorenstein global
   dimension. If $\mathcal{SGF}(R)$ is projectively resolving then $\mathcal{SGP}(R)$ is projectively
   resolving, with equivalence if $R$ is perfect.
   \end{prop}
   \proof Assume that $R$ is a ring with finite Gorenstein global
   dimension such that $\mathcal{SGF}(R)$ is projectively
   resolving. To prove that $\mathcal{SGP}(R)$ is resolving, we should, from
   Lemma \ref{lem SGP}, to prove that every $G$-projective $R$-module is
   $SG$-projective. So, let $M$ be a $G$-projective module. By \cite[Theorem 2.20]{Holm}, $Ext_R(M,Q)=$ for every
   projective module $Q$. So, to prove that $M$ is $SG$-projective it suffices to find a short exact sequence
   of $R$-modules $0  \longrightarrow M \longrightarrow P \longrightarrow M \longrightarrow 0$ where $P$ is projective
   (by \cite[Proposition 2.9]{Bennis and Mahdou1}). From
   \cite[Proposition 3.4 and Theorem 3.24]{Holm}, $M$ is also
   $G$-flat (since $G-gldim(R)$ is finite). Then, $M$ is $SG$-flat
   since $\mathcal{SGF}(R)$ is projectively resolving (by Lemma \ref{lem
   SGF}). Thus, from \cite[Proposition 3.6]{Bennis and Mahdou1},
   there exists a short exact sequence $0\longrightarrow N \longrightarrow F \longrightarrow M
   \longrightarrow 0$ where $F$ is flat. So, from \cite[Corollary 2.10]{Bennis and
   Mahdou2}, $pd_R(F)$ is finite. On the other hand, by \cite[Theorem
   3.7]{Holm}, $\mathcal{GF}(R)$ is projectively resolving and then, from the
   short exact sequence above, $F$ is $G$-projective since $M$ is
   $G$-projective. Therefore, from \cite[Proposition 2.27]{Holm},
   $F$ is projective. Consequently, we have the desired short
   exact sequence.\\
   Now, assume that $R$ is a perfect ring with finite Gorenstein
   global dimension and such $\mathcal{SGP}(R)$ is projectively resolving.
   we claim the projectively resolving of the class
   $\mathcal{SGF}(R)$. So, to prove this, it suffices to prove
   that every $G$-flat module $M$ is $SG$-flat. By \cite[Theorem
   3.14]{Holm}, $Tor_R(M,I)=0$ for every injective module $I$. So,
   from \cite[Proposition 3.6]{Bennis and Mahdou1}, it stays to
   prove the existence of a short exact sequence $0
   \longrightarrow M \longrightarrow F \longrightarrow M
   \longrightarrow 0$ where $F$ is flat. By \cite[Theorem 3.5]{Bennis and
   Mahdou1} $M$ is a direct summand of an $SG$-flat module $N$.
   For such module, by \cite[Proposition 3.6]{Bennis and Mahdou1}, there exists a short exact sequence $(\ast) \quad 0
   \longrightarrow N \longrightarrow F \longrightarrow N
   \longrightarrow 0$ where $F$ is flat (then projective since $R$
   is perfect). Now, let $P$ be a  projective module.
   We have $id_R(P) \leq n$ where $n=G-gldim(R)$ (by \cite[Corollary 2.10]{Bennis and
   Mahdou2}). Then, $Ext_R^{n+1}(N,P)=0$ and so, from the short exact
   sequence $(\ast)$ we deduce that $Ext_R(N,P)=0$. So, from
   \cite[Proposition 2.9]{Bennis and Mahdou1}, $N$ is
   $SG$-projective. Consequently, $M$ is $G$-projective (since it is
   a direct summand of $N$ and by \cite[Theorem 2.7]{Bennis and
   Mahdou1}). Then, $M$ is $SG$-projective since $\mathcal{SGP}(R)$ is
   projectively resolving and by Lemma \ref{lem SGP}. Hence, by \cite[Proposition 2.9]{Bennis and Mahdou1}, there
   exists a short exact sequence $0\longrightarrow M
   \longrightarrow P  \longrightarrow M \longrightarrow 0$ where
   $P$ is projective (then  flat), and this is the desired short
   exact sequence.$\blacksquare$

\end{section}




\end{document}